\documentclass[a4paper]{amsart}  

\usepackage{amssymb} \usepackage{amscd} \usepackage{times}

\def\zbb{\mathbb{Z}}  
  
  \def\phi{\varphi}
 \def\p1{{\mathbb{P}^1_\zbb}}

\begin{document}

\title{ Harnack inequalities for Yamabe type equations}

\author{Samy Skander Bahoura}

\address{Department of Mathematics, Patras University, 26500 Patras , Greece }
              
\email{samybahoura@yahoo.fr, bahoura@ccr.jussieu.fr} 

\date{}

\maketitle

\begin{abstract}

We give some a priori estimates of type $ \sup \times \inf $ on Riemannian manifolds for Yamabe  and prescribed curvature type equations. An application of those results is the uniqueness result for $ \Delta u+\epsilon u=u^{N-1} $ with $ \epsilon $ small enough.

\end{abstract}

\bigskip

\bigskip

\begin{center}  INTRODUCTION AND RESULTS.

\end{center}

\bigskip

\bigskip

We are on Riemannian manifold $ (M,g) $ of dimension $ n \geq 3 $. In this paper we denote $ \Delta = -\nabla^j(\nabla_j) $ the geometric laplacian and $ N=\dfrac{2n}{n-2} $.

\smallskip

The scalar curvature equation is:

$$ \dfrac{4(n-1)}{n-2}\Delta u+R_gu=V u^{N-1},\,\, u >0 . $$

Where $ R_g $ is the scalar curvature and $ V $ is a function (prescribed scalar curvature).

\bigskip

When we suppose $ V \equiv 1 $, the previous equation is the Yamabe equation.

\bigskip

Here we study some properties of Yamabe and prescribed scalar curvature equations. The existence result for the Yamabe equation on compact Riemannian manifolds was proved by T. Aubin and R. Schoen ( see for example [Au]).

\bigskip

First, we suppose the manifold $ (M,g) $ compact. We have:

\bigskip

{\it Theorem 1. For all $ a,b,m>0 $, there exist a positive constant $ C=C(a,b,m,M,g) $ such that for every $ \epsilon >0 $, for every smooth function $ V $ such that  $  a\leq V_{\epsilon}(x) \leq b $ and every positive solution $ u_{\epsilon} $ of:

$$ \Delta u_{\epsilon}+\epsilon u_{\epsilon}=V_{\epsilon} {u_{\epsilon}}^{N-1} $$

with $ \max_M u_{\epsilon} \geq m $, we have:

$$ \epsilon \max_M u_{\epsilon} \min_M u_{\epsilon} \geq C. $$}

Now, we  consider a Riemannian manifold $ (M,g) $ of dimension $ n\geq 3 $ ( not necessarily compact) and we work with Yamabe type equation,

$$ \Delta u-\lambda u=n(n-2)u^{N-1}. $$

We look for a priori bounds for solutions of the previous equation.

\bigskip

{\it Theorem 2. If $ 0 < m \leq \lambda+R_g \leq 1/m $ then for every compact $ K $ of $ M $, there exist a positive constant $ c=c(K,M,m,n,g) $ such that:

$$ \sup_K u \times \inf_M u \leq c. $$}

Note that there is lot of estimates of those type for prescribed scalar curvature on open set $ \Omega $ of $ {\mathbb R}^n $, see ([B],[B-M], [B-L-S], [C-L 1], [C-L 2], [L 1], [L 2], and [S]). 

\bigskip

In dimension 2 Brezis, Li and Shafrir [B-L-S], have proved that $ \sup + \inf $ is bounded from above when we suppose the prescribed curvature uniformly lipschitzian. In [S], Shafrir got a result of type $ \sup +C\inf $, with $ L^{\infty} $ assumption on prescribed curvature.

\bigskip

In dimensions $ n\geq 3 $, we can find many results with different assumptions on prescribed curvature, see [B], [L 2], [C-L 2].

\bigskip

Note that an important estimates was proved for Yamabe equation about the product $ \sup \times \inf $, in dimensions 3,4 by Li and Zhang [L-Z].

\bigskip

In our work we have no assumption on energy. There is an important work if we suppose the energy bounded, see for example [D-H-R].

\bigskip

\underbar {\bf Application:}

\bigskip

We assume that $ M $ is compact and $ 1/m \geq R_g \geq m >0 $ on $ M $. For small values of $ \lambda $ we can have some upper bounds for the product $ \sup \times \inf $ for the following equation:

$$ \Delta u_{\epsilon}+ \epsilon u_{\epsilon}=n(n-2)u_{\epsilon}^{N-1}. $$

{\it Theorem 3. If $ \epsilon \to 0 $, then,

$$ \sup_M u_{\epsilon} \times \inf_M u_{\epsilon} \leq c(n,m,M,g). $$}

A consequence of Theorems 1 and 3 is the following corollary:

\bigskip

{\it Corollary. Any sequence $ u_i>0 $ solutions of the following equation:

$$ \Delta u_i+\epsilon_i u_i=n(n-2){u_i}^{N-1} , $$ 

 converge uniformly  to 0 on $ M $ when $ \epsilon_i $ tends to $ 0 $.}

\bigskip

We have:

\bigskip

{\it Theorem 4. On compact Riemannian manifold $ (M,g) $ with $ R_g >0 $ every-where, the sequence $ u_i >0 $ solutions of the previous equation is such that for $ i $ large, $ u_i\equiv \left [ \dfrac{\epsilon_i}{n(n-2)} \right ]^{(n-2)/4} $.}

\bigskip

Note that the previous result assert that $ \left [ \dfrac{\epsilon_i}{n(n-2)} \right ]^{(n-2)/4} $ is the only solution of the previous equation for $ \epsilon_i $ small.

\bigskip

We remark an important result in [B-V,V]; they have a same consequence than in theorem 4 with assumption on Ricci curvature ( $ Ric \geq \epsilon_i c g $, with $ c>0 $). Here we give a condition on scalar curvature to obtain an uniqueness result.

\newpage

\underbar {\bf Proof of theorem 1:}

\bigskip

We need two lemmata and one proposition. We are going to prove some estimates for the Green function $ G_{\epsilon} $ of the operator $ \Delta + \epsilon $.

\bigskip

\underbar {\bf Lemma 1.}

\bigskip

 For each point $ x \in M $ there exist $ \epsilon_0 >0 $ and $ C(x,M,g)>0 $ such that for every $ z \in B(x_0, \epsilon_0) $, and every $ \mu \leq \epsilon_0 $, every $ a,b \in \partial B(z, \mu) $, there exist a curve $ \gamma_{a,b} $ of classe $ C^1 $ linking $ a $ to $ b $ which included in $ \partial B(z,\mu) $. The length of this curve is  $ l(\gamma_{a,b}) \leq C(x,M,g)\mu $.

\bigskip

\underbar {\bf Proof:}

\bigskip

Let $ x \in M $, we consider a chart $ (\Omega, \phi) $ around $ x $.

\smallskip

We take exponential map on the compact manifold $ M $. According to T. Aubin and E. Hebey see [Au] and [He], there exist $ \epsilon >0 $ such that $ \exp_x $ is $ C^{\infty} $ function of $ B(x,\epsilon) \times B(0,\epsilon) $ into $ M $ and for all $ z\in B(x,\epsilon) $, $ \exp_z $ is a diffeomorphism from $ B(0,\epsilon) $ to $ B(z, \epsilon) $ with $ \exp_z[\partial B(0,\mu)]=\partial B(z,\mu) \subset M $ for $ \mu \leq \epsilon/2 $. If we take two points $ a,b $ of $ \partial B(z,\mu) $ ($\mu \leq  \epsilon/2$ ), then $ a'=\exp_z^{-1}(a), b'=\exp_z^{-1}(b) $ are two points of $ \partial B(0,\mu) \subset {\mathbb R}^n $. On this sphere of center 0 and radius $ \mu $, we can link $ a' $ to $ b' $ by a great circle arc whose length is $ \leq 2\pi\mu $. Then, there exist a curve of class $ C^1 $ $ \delta_{a',b'} $ in $ \partial B(0,\mu) \subset {\mathbb R}^n $ such that $ l(\delta_{a',b'})\leq 2\pi\mu $. Now we consider the curve $ \gamma_{a,b}=\exp_z(\delta_{a',b'}) $, this curve of class $ C^1 $, link $ a $ to $ b $ and it is included in $ \partial B(z,\mu) \subset M $. The length of $ \gamma_{a,b} $ is giving by the following formula :

$$ l(\gamma_{a,b})=\int_0^1\sqrt {g_{ij}[\gamma_{a,b}(s)](\dfrac{d\gamma_{a,b}}{dt})^i(s)(\dfrac{d \gamma_{a,b}}{dt})^j(s)}ds. $$

where $ g_{ij} $ is the local expression of the metric $ g $ in the chart $ (\Omega, \phi) $.

\bigskip

We know that there exist a constant $ C=C(x,M,g)>1 $ such that:

$$ \dfrac{1}{C}||X||_{{\mathbb R}^n} \leq g_{ij}(z) X^iX^j \leq C||X||_{{\mathbb R}^n} \,\, {\rm for \,\, all } \,\, z\in B(x,\epsilon/2)\,\,{\rm and}\,\, {\rm all } \,\, X \in {\mathbb R}^n . $$

 We have,

$$ l(\gamma_{a,b})=\int_0^1 \sqrt { g_{ij}[\exp_z[\delta_{a',b'}(s)]] (\dfrac{ d [\exp_z[\delta_{a',b'}]]}{dt})^i(s) (\dfrac{ d[\exp_z[\delta_{a',b'}]]}{dt})^j(s) }ds, $$

$$ l(\gamma_{a,b}) \leq C
\int_0^1 ||d\exp_z(\dfrac{d\delta_{a',b'}}{dt})(s) ||_{{\mathbb R}^n}ds, $$

and $ u:(z,v)\to \exp_z(v) $ is $ C^{\infty} $ on $ B(x,\epsilon) \times B(0,\epsilon) $, 

\bigskip

but,

$$ || du_{z,v} || = ||d \exp_z(v)|| \leq C'(x,M,g) \,\,\, \forall \,\, (z,v)\in B(x,\epsilon/2)\times B(0,\epsilon/2) {\rm ( in\,\, the \,\, sense \,\, of \,\, linear \,\, form)}, $$

Finaly, 

$$ l(\gamma_{a,b})\leq  \tilde C(x,M,g) \int_0^1 ||\dfrac{ d\delta_{a',b'}}{dt}(s) ||_{{\mathbb R}^n} ds=\tilde C(x,M,g)l(\delta_{a',b'})\leq 2\pi \tilde C(x,M,g) \mu . $$

We need to estimate the singularities of Green functions. Set $ G_i=G_{\epsilon_i} $.

\bigskip

\underbar {\bf Lemma 2.} 

\bigskip

The function $ G_i $ satisfies:

$$ G_i(x,y)\leq \dfrac{C(M,g)}{\epsilon_i [d_g(x,y)]^{n-2}}. $$

where $ C(M,g)>0 $ and $ d_g $ is the distance on $ M $ for the metric $ g $.

\bigskip

\underbar {\bf Proof:}

\smallskip
 
According to the Appendix of [D-H-R] (see also [Au]), we can write the function $ G_i $ :

$$ G_i(x,y)=H(x,y)+\Sigma_{j=1}^k \Gamma_{i,k}(x,y)+u_{i,k+1}(x,y), $$

with, $ k=[n/2] $ and $ u_{i,k+1} $ is solution of $ \Delta u_{i,k+1} + \epsilon_i u_{i,k+1}= \Gamma_{i,k+1} $.

\smallskip
 
According to Giraud (see [Au] and [D-H-R]), we have:

\smallskip

 i) $ 0\leq H(x,y)\leq \dfrac{C_0(M,g)}{[d_g(x,y)]^{n-2}} , $

\smallskip

ii) $ |\Gamma_{i,j}(x,y)| \leq \dfrac{C_j(M,g)}{[d_g(x,y)]^{n-2}} , j=1,\ldots, k $ and,

\smallskip

iii) $ \Gamma_{i,k+1}(x,y) \leq C_{k+1}(M,g) $ and continuous on $ M \times M $ .

\bigskip

We write $ u_{i,k+1} $ by using the Green function  $ G_i $, we obtain with iii):

$$ u_{i,k+1}(x,y)=\int_M G_i(x,y)\Gamma_{i,k+1}(x,y)dV_g(y) \leq C_{k+1}(M,g)\int_M G_i(c,y)dV_g(y)=\dfrac{C_{k+1}(M,g)}{\epsilon_i}. $$

If we combine the last inequality and i) et ii), we obtain the result of the lemma.

\bigskip

We have to estimate the Green function from below.

\underbar {\bf Proposition.} 

\bigskip

Consider two sequences of points of $ M $, $ ( x_i) $ et $ (y_i) $ such that $ x_i\not = y_i $ for all $ i $ and $ x_i \to x $, $ y_i \to y $. Then, there exist a positive constant $ C $ depending on $ x,y, M $ and $ g $, and a subsequence $ (i_j ) $ such that:
 
$$ G_{i_j}(x_{i_j},y_{i_j}) \geq \dfrac{C}{\epsilon_{i_j}} \qquad  \forall \,\, j . $$

\underbar {\bf Proof:}

\bigskip

We know that $ G_i(x_i,.) $ is $ {C}^{\infty}(M-{x_i}) $ and satisfies the following equation:

$$ \Delta G_i(x_i,.)+\epsilon_i G_i(x_i,.)=0, \,\,\, {\rm in} \,\, M-{x_i} . $$

{\underbar {Case 1:} }$ y=x $.

Let $ R_i=\dfrac{1}{2} d_g(x_i,y_i)>0 $ and $ \Omega_i=M-B(x_i,R_i) $, according to maximum principle, the function $ G_i(x_i,.) $ has its maximum on the boundary of $ \Omega_i $. Then;

$$ \max_{\Omega_i} G_i(x_i,z)=G_i(x_i,z_i), \,\, d(x_i,z_i)=R_i . $$

Let $ t_i $ be a point of $ M $ such that $ d_g[y_i,B(x_i,R_i)]=d(t_i,y_i) $. We have $ t_i \in \partial B(x_i,R_i) $  then $ d(x_i,t_i)=R_i $. Because the manifold $ M $ is compact, we can find a minimizing curve $ L_i $ between $ y_i $ and $ t_i $. Let $ \delta_i $ a curve in $ \partial B(x_i,R_i) $ with minimal length linking $ t_i $ to $ z_i $. We can choose it like in lemma 1. Then $ l(\delta_i) \leq c(x,M,g) R_i $ and if we note $ \bar \delta_i=\delta_i \cup L_i $, we have $ l(\bar \delta_i)=l(\delta_i)+l(L_i) \leq R_i[1+c(x,M,g)] $. The curve $ \bar \delta_i $ link $ z_i $ to $ y_i $, and it is included in $ \Omega_i $. Let $ r_i=\dfrac{1}{5} R_i $. We cover the curve $ \bar \delta_i $ by balls of radii $ r_i $, if we consider $ N_i $ the minimal number  of those balls, then we have $ N_i r_i \leq [c(x,M,g)+1] R_i $, and  $ N_i \leq 5[c(x,M,g)+1] $.

\bigskip

If we work on open set of one chart $ \Omega $ centered in $ x $, with a small ball around $ x_i $ removed, $ \tilde \Omega_i=\Omega-{B(x_i,\dfrac{1}{100} R_i)} $ then, we can apply the theorem 8.20 of [GT] (Harnack inequality)  in each ball of the finite covering of $ \bar \delta_i $ defined previous. In this Harnack inequality the constant which depends on the radius is explicit and equal to $ C_0(n)^{(\Lambda/\lambda) + \nu R_i)}  $ but here $ R_i \to 0 $, and the constant do not depend on the radius. We obtain:

$$ \sup_{B(z_i,r_i)} G_i(x_i,z) \leq C(x,M,g) \inf_{B(y_i,r_i)} G_i(x_i,z) . $$

Then, $ G_i(x_i,z)\leq G_i(x_i,z_i) \leq C(x,M,g) G_i(x_i,y_i) $ for all $ z \in \Omega_i $.

\bigskip

Now we write:

$$\dfrac{1}{\epsilon_i}=\int_M G_i(x_i,z)dV_g(z)=\int_{\Omega_i} G_i(x_i,z)dV_g(z)+\int_{B(x_i,2R_i)} G_i(x_i,z) dV_g(z), $$

but,

$$ \int_{\Omega_i} G_i(x_i,z)dV_g(z) \leq |\Omega_i| \sup_{\Omega_i} G_i(x_i,z) \leq |\Omega_i| C(x,M,g) G_i(x_i,y_i), $$

we take $ A_i= \int_{B(x_i,2R_i)} G_i(x_i,z)dV_g(z) $, we have,

$$ A_i=\int_{B(0,2R_i)} G_i[x_i,\exp_{x_i}(v)]{ \sqrt {|g|}} du = \int_0^{2R_i} \int_{{\mathbb S}_{n-1}}t^{n-1}{\sqrt {|g|}} G_i[x_i,\exp_{x_i}(t\theta)] dtd\theta, $$

if we use the lemma 8 in  Hebey-Vaugon (see [H-V]), we obtain $ {\sqrt {|g|}} \leq c(M,g) $. But $ R_i \to 0 $, then $ d_g[x_i,\exp_{x_i}(t\theta)]=t $ (the geodesic are minimizing). We use the lemma 2 and we find:

$$ \int_{B(x_i,2R_i)} G_i(x_i,z)dV_g(z) \leq \dfrac{c'(M,g)(R_i)^2}{\epsilon_i}.$$

Finaly:

$$ G_i(x_i,y_i) \leq \dfrac{1-c'(M,g)(R_i)^2}{|\Omega_i|C(x,M,g)\epsilon_i} \geq \dfrac{C'(x,M,g)}{\epsilon_i}. $$

\bigskip

{\underbar {Case 2}}: $ x \not = y $.

\bigskip

We write,

$$ \dfrac{1}{\epsilon}=\int_M G_i(x_i,z)dV_g(z)=\int_{M-B(x_i,\delta)} G_i(x_i,z)dV_g(z)+\int_{B(x_i,\delta)} G_i(x_i,z)dV_g(z). $$

We take $ 0 < \delta \leq \dfrac{inj_g(M)}{2} $, where $ inj_g(M) $ is the injectivity radius of the compact manifold $ M $. We use the exponential map and we have:

$$ \int_{B(x_i,\delta)} G_i(x_i,z)dV_g(z)=\int_{B(0,\delta) }G_i[x_i,\exp_{x_i}(v)] \sqrt {|g|} dv=\int_0^{\delta} t^{n-1} \int_{{\mathbb S}_{n-1}} G_i[x_i, \exp_{x_i}(t\theta)] \sqrt |g|dtd\theta, $$

If we use the lemma 8 in Hebey-Vaugon (see [H-V]), we obtain $ |g|\leq c(M,g) $. Using the fact $ t\to \exp_{x_i}(t\theta) $ is minimizing for $ t\leq \delta < inj_g(M) $ and  the lemma 2, we obtain:

$$ \int_{B(x_i,\delta)} G_i(x_i,z)dV_g(z) \leq \dfrac{ C'(M,g) \delta^2}{\epsilon_i} . $$

Then,

$$ \int_{M-B(x_i,\delta)}G_i(x_i,z)dV_g(z) \geq \dfrac{1-C'(M,g)\delta^2}{\epsilon_i}, $$

we can choose $ 0 < \delta < \dfrac{1}{\sqrt {C'(M,g)}} $.

\bigskip

Between $ x $ and $ y $, we work like in the first case. We take $ 0 < \delta < \dfrac{ d_g(x,y)}{2} $, for each $ i $ and we consider the maximum of $ G_i(x_i,.) $ in $ \Omega_i=M-B(x_i,\delta) $. By maximum principle $ G_i(x_i,z_i)=\max_{\Omega_i} G_i(x_,z)=\max_{\partial B(x_i,\delta)} G_i(x_i,z) $. After passing to a subsequence, we can suppose that $ z_i \to z $.

\bigskip

We have $ 0 < \delta=d(x_i,z_i)\to d(x,z) $. We choose $ \delta >0 $ such that the ball of center $ x $ and radius $ 2\delta $ is included in open chart centred in $ x $. (we can choose the exponential map in $ x $ and use the lemma 1). 

\bigskip

Let $ t $ be the point of $ B(x,\delta) $ such that $ d(y,t)=d[y,B(x,\delta)] $, $ t $ depend on $ x $ and $ y $. We consider a minimizing curve $ L_1 $ between $ t $ and $ y $. The manifold is compact and  $ \delta << inj_g(M) $, then, in each point $ u $ of $ L_1 $, $ [ B(u,\delta/2), exp_u ] $ is a local chart. We cover the curve $ L_1 $ by a finite number of balls of radii $ \delta/10 $. We apply the Harnack inequality between those balls for the functions $ G_i(x_i,.) $. We infer that:

$$ \sup_{B(t,\delta/10)} G_i(x_i,s) \leq C(x,y,M,g) \inf_{B(y_i,\delta/10)} G_i(x_i,s) \leq C(x,y,M,g) G_i(x_i,y_i) . $$

Now we want to know what happens between the balls $ B(t,\delta/10) $ and $ B(z,\delta/10) $. The ball $ B(x,2\delta) $ is open chart set centered in $ x $. We choose a curve $ L_2 $ between $ z $ and $ t $ like in the first case. This curve must stay in $ \partial B(x,\delta) $ and its length $ l \leq C_1(x,y,M,g) \delta $, then, we can have a covring of this curve by a minimal number $ N $ of balls od radii $ \delta/10 $, in fact $ N\leq C_2(x,y,M,g) $ (like in the first case). Those balls are included in the open chart set centered in $ x $ which we choose as in the begining. Then, the operator $ \Delta+\epsilon_i $ has those coefficients depending only on the open chart set centred in $ x $  and not depending on  $ z $, we can apply the Harnack inequality (theorem 8.20 of [GT]) in this open set without $ B(x,\delta/100) $, for the functions $ G_i(x_i,.) $. Finaly, we obtain the same conclusion than in the case 1, there exist $ C=C(x,y,M,g)>0 $ such:

$$ G_i(x_i,s)\leq C G_i(x_i,y_i) \,\,\, \forall \,\, s\in \Omega_i=M-B(x_i,\delta) \,\, \forall \,\, i \geq i_0.$$

The rest of the proof is the same than in the case 1.

\bigskip

\underbar {\it Proof of Theorem 1.}

\bigskip

We write $ u_i $ by using the Green function $ G_i $, then:

$$ \min_M u_i=u_i(x_i)=\int_M G_i(x_i,y) V_i(y){u_i(y)}^{N-1}dV_g(y), $$

then,

$$ \sup_M u_i \times \inf_M u_i \geq \int_M G_i(x_i,y) V_i(y) {u_i(y)}^N dV_g(y) \geq a\min_M G_i(x_i,.) \int_M {u_i(y)}^N dV_g(y) .$$

Let $ G_i(x_i,y_i)=\min_M G_i(x_i,.) $, after passing to a subsequence, we can assume that $ x_i \to x $ and $ y_i \to y $. By using the previous proposition, we can suppose that  there exist a positive constant $ c=c(x,y,M,g) $ such that:

$$ G_i(x_i,y_i) \geq \dfrac{c}{\epsilon_i}. $$

Then,

$$ \int_M [u_i(y)]^N dV_g(y) \leq \epsilon_i \sup_M u_i \inf_M u_i. $$

We know argue by contradiction and assume that $ \epsilon_i \sup_M u_i \times \inf_M u_i $ tends to 0. We know ( see a previous paper when we use the Moser iterate scheme, see [B1]), that ( after passing to a subsequence) for $ q $ large:

$$  ||u_i||_{L^q(M)} \to 0. $$

Assume that $ G $ the Green function of the laplacian, we can write:

$$ u_i(x)=\dfrac{1}{Vol(M)} \int_M u_i + \int_M G(x,y)[V_i(y) {u_i(y)}^{N-1}-\epsilon_iu_i(y)] dV_g(y), $$

and if we use Holder inequality, we obtain:

$$ \sup_M u_i \to 0 .$$

But, this is a contradiction with $ \sup_M u_i \geq m>0 . $

\newpage

\underbar {\bf Proof of the theorems 2,3,4.}

\bigskip

{\underbar {\bf Part I: The metric in polar coordinates.}}

\bigskip

Let $ (M,g) $ a Riemannian manifold. We note $ g_{x,ij} $ the local expression of the metric $ g $  in the exponential map centred in $ x $.

We are concerning by the polar coordinates expression of the metric. Using Gauss lemma, we can write:

$$ g=ds^2=dt^2+g_{ij}^k(r,\theta)d\theta^id\theta^j=dt^2+r^2{\tilde g}_{ij}^k(r,\theta)d\theta^id\theta^j=g_{x,ij}dx^idx^j, $$

in a polar chart with origin $ x $", $ ]0,\epsilon_0[\times U^k $, with $ ( U^k, \psi) $ a chart of $ {\mathbb S}_{n-1} $. We can write the element volume:

$$ dV_g=r^{n-1}\sqrt {|{\tilde g}^k|}dr d \theta^1 \ldots d \theta^{n-1} = \sqrt {[det(g_{x,ij})]}dx^1 \ldots dx^n, $$

then,

$$ dV_g=r^{n-1} \sqrt { [det(g_{x,ij})]}[\exp_x(r\theta)]\alpha^k(\theta)dr d\theta^1 \ldots d \theta^{n-1} , $$

where, $ \alpha^k $ is such that, $ d\sigma_{{\mathbb S}_{n-1}}=\alpha^k(\theta) d \theta^1 \ldots d \theta^{n-1} . $ (Riemannian volume element of the sphere in the chart $ (U^k,\psi) $ ).

\bigskip

Then,

$$ \sqrt { |{\tilde g}^k|}=\alpha^k(\theta) \sqrt {[det(g_{x,ij})]}. $$

Clearly, we have the following proposition:

\bigskip

\underbar {\bf Proposition 1:} Let $ x_0 \in M $, there exist $ \epsilon_1>0 $ and if we reduce $ U^k $, we have:

$$ |\partial_r{\tilde g}_{ij}^k(x,r,\theta)|+|\partial_r\partial_{\theta^m}{\tilde g}_{ij}^k(x,r, \theta)| \leq C r,\,\, \forall \,\, x\in B(x_0,\epsilon_1) \,\, \forall \,\, r\in [0,\epsilon_1], \,\, \forall \,\, \theta \in U^k.$$

and,

$$ |\partial_r|{\tilde g}^k|(x,r,\theta)|+\partial_r \partial_{\theta^m} |{\tilde g}^k|(x,r,\theta)\leq C r,\,\, \forall \,\, x\in B(x_0,\epsilon_1) \,\, \forall \,\, r\in [0,\epsilon_1], \,\, \forall \,\, \theta \in U^k. $$

\underbar {\bf Remark:} 

\bigskip

$ \partial_r [ \log \sqrt { |{\tilde g}^k |}] $ is a local function of $ \theta $, and the restriction of the global function on the sphere $ {\mathbb S}_{n-1} $, $ \partial_r [\log \sqrt { det(g_{x, ij})}] $. We will note, $ J(x,r,\theta)=\sqrt { det(g_{x, ij})} $.

\bigskip

{\underbar { \bf Part II: The laplacian in polar coordinates}}

\bigskip

Let's write the laplacian in $ [0,\epsilon_1]\times U^k $,

$$ -\Delta = \partial_{rr}+\dfrac{n-1}{r}\partial_r+ \partial_r [\log \sqrt { |{\tilde g^k|}] }\partial_r+\dfrac{1}{r^2 \sqrt {|{\tilde g}^k|}}\partial_{\theta^i}({\tilde g}^{\theta^i \theta^j}\sqrt { |{\tilde g}^k|}\partial_{\theta^j}) . $$

We have,

$$ -\Delta = \partial_{rr}+\dfrac{n-1}{r}\partial_r+ \partial_r \log J(x,r,\theta)\partial_r+ \dfrac{1}{r^2 \sqrt {|{\tilde g}^k|}}\partial_{\theta^i}({\tilde g}^{\theta^i \theta^j}\sqrt { |{\tilde g}^k|}\partial_{\theta^j}) . $$

We write the laplacian ( radial and angular decomposition),

$$ -\Delta = \partial_{rr}+\dfrac{n-1}{r} \partial_r+\partial_r [\log J(x,r,\theta)] \partial_r-\Delta_{{S}_r(x)}, $$

where $ \Delta_{ S_r(x)} $ is the laplacian on the sphere $ {S}_r(x) $. 

\bigskip

We set $ L_{\theta}(x,r)(...)=r^2\Delta_{ S_r(x)}(...)[\exp_x(r\theta)] $, clearly, this operator is a laplacian on $ {\mathbb S}_{n-1} $ for particular metric. We write,

$$ L_{\theta}(x,r)=\Delta_{g_{x,r, {}_{{\mathbb S}_{n-1}}}}, $$

and,

$$ \Delta = \partial_{rr}+\dfrac{n-1}{r} \partial_r+\partial_r [ J(x,r,\theta)] \partial_r - \dfrac{1}{r^2} L_{\theta}(x,r) . $$

If, $ u $ is function on $ M $, then, $ \bar u(r,\theta)=u[\exp_x(r\theta)] $ is the corresponding function in polar coordinates centred in $ x $. We have,

$$ -\Delta u =\partial_{rr} \bar u+\dfrac{n-1}{r} \partial_r \bar u+\partial_r [ J(x,r,\theta)] \partial_r \bar u-\dfrac{1}{r^2}L_{\theta}(x,r)\bar u . $$

{\underbar { \bf Part III: "Blow-up" and "Moving-plane" methods }} 

\bigskip

\underbar {\bf The "blow-up" technic}

\bigskip

Let, $ (u_i)_i $ a sequence of functions on $ M $ such that,

$$ \Delta u_i -\lambda u_i =n(n-2){u_i}^{N-1}, \,\, u_i>0,\,\, N=\dfrac{2n}{n-2}, \qquad (E)  $$

We argue by contradiction and we suppose that $ \sup \times \inf $ is not bounded.

\smallskip

We assume that:

\bigskip

$ \forall \,\, c,R >0 \,\, \exists \,\, u_{c,R} $ solution of $ (E) $ such that:

$$ R^{n-2} \sup_{B(x_0,R)} u_{c,R} \times \inf_M u_{c,R} \geq c. \qquad (H) $$

\bigskip

\underbar {\bf Proposition 2:} 

\bigskip

There exist a sequence of points $ (y_i)_i $, $ y_i \to x_0 $ and two sequences of positive real number $ (l_i)_i, (L_i)_i $, $ l_i \to 0 $, $ L_i \to +\infty $, such that if we consider $ v_i(y)=\dfrac{u_i[\exp_{y_i}(y)]}{u_i(y_i)} $, we have:

$$ i) \qquad 0 < v_i(y) \leq  \beta_i \leq 2^{(n-2)/2}, \,\, \beta_i \to 1. $$

$$ ii) \qquad v_i(y)  \to \left ( \dfrac{1}{1+{|y|^2}} \right )^{(n-2)/2}, \,\, {\rm uniformly \,\, on\,\, every \,\, compact \,\, set \,\, of } \,\, {\mathbb R}^n . $$

$$ iii) \qquad l_i^{(n-2)/2} [u_i(y_i)] \times \inf_M u_i \to +\infty $$

\underbar {\bf Proof:}

We use the hypothesis $ (H) $. We can take two sequences $ R_i>0, R_i \to 0 $ and $ c_i \to +\infty $, such that,

$$ {R_i}^{(n-2)} \sup_{B(x_0,R_i)} u_i \times \inf_M u_i \geq c_i \to +\infty. $$

Let, $ x_i \in  { B(x_0,R_i)} $, such that $ \sup_{B(x_0,R_i)} u_i=u_i(x_i) $ and $ s_i(x)=[R_i-d(x,x_i)]^{(n-2)/2} u_i(x), x\in B(x_i, R_i) $. Then, $ x_i \to x_0 $.

\bigskip

We have, 

$$ \max_{B(x_i,R_i)} s_i(x)=s_i(y_i) \geq s_i(x_i)={R_i}^{(n-2)/2} u_i(x_i)\geq \sqrt {c_i}  \to + \infty. $$ 

\bigskip

Set :

$$ l_i=R_i-d(y_i,x_i),\,\, \bar u_i(y)= u_i [\exp_{y_i}(y)],\,\,  v_i(z)=\dfrac{u_i [ \exp_{y_i}\left ( z/[u_i(y_i)]^{2/(n-2)} \right )] } {u_i(y_i)}. $$

Clearly, $ y_i \to x_0 $. We obtain:

$$ L_i= \dfrac{l_i}{(c_i)^{1/2(n-2)}} [u_i(y_i)]^{2/(n-2)}=\dfrac{[s_i(y_i)]^{2/(n-2)}}{c_i^{1/2(n-2)}}\geq \dfrac{c_i^{1/(n-2)}}{c_i^{1/2(n-2)}}=c_i^{1/2(n-2)}\to +\infty. $$

\bigskip

If $ |z|\leq L_i $, then $ y=\exp_{y_i}[z/ [u_i(y_i)]^{2/(n-2)}] \in B(y_i,\delta_i l_i) $ with $ \delta_i=\dfrac{1}{(c_i)^{1/2(n-2)}} $ and $ d(y,y_i) < R_i-d(y_i,x_i) $, thus, $ d(y, x_i) < R_i $ and, $ s_i(y)\leq s_i(y_i) $, we can write,

$$ u_i(y) [R_i-d(y,y_i)]^{(n-2)/2} \leq u_i(y_i) (l_i)^{(n-2)/2}. $$

But, $ d(y,y_i) \leq \delta_i l_i $, $ R_i >l_i$ and $ R_i-d(y, y_i) \geq R_i-\delta_i l_i>l_i-\delta_i l_i=l_i(1-\delta_i) $, we obtain,

$$ 0 < v_i(z)=\dfrac{u_i(y)}{u_i(y_i)} \leq \left [ \dfrac{l_i}{l_i(1-\delta_i)} \right ]^{(n-2)/2}\leq 2^{(n-2)/2} . $$

We set, $ \beta_i=\left ( \dfrac{1}{1-\delta_i} \right )^{(n-2)/2} $, clearly $ \beta_i \to 1 $.

\bigskip

The function $ v_i $ is solution of:

$$ -g^{jk}[\exp_{y_i}(y)]\partial_{jk} v_i-\partial_k \left [ g^{jk}\sqrt { |g| } \right ][\exp_{y_i}(y)]\partial_j v_i+ \dfrac{ R_g[\exp_{y_i}(y)]}{[u_i(y_i)]^{4/(n-2)}} v_i=n(n-2){v_i}^{N-1}, $$

By elliptic estimates and Ascoli, Ladyzenskaya theorems, $ ( v_i)_i $ converge uniformely on each compact to the function $ v $ solution on $ {\mathbb R}^n $ of, 

$$ \Delta v=n(n-2)v^{N-1}, \,\, v(0)=1,\,\, 0 \leq v\leq 1\leq 2^{(n-2)/2}, $$

By using maximum principle, we have $ v>0 $ on $ {\mathbb R}^n $, the result of Caffarelli-Gidas-Spruck ( see [C-G-S]) give, $ v(y)=\left ( \dfrac{1}{1+{|y|^2}} \right )^{(n-2)/2} $. We have the same properties for $ v_i $ in the previous paper [B2].

\bigskip

\underbar {\bf Polar coordinates and "moving-plane" method}

\bigskip

Let, 

$$ w_i(t,\theta)=e^{(n-2)/2}\bar u_i(e^t,\theta) = e^{(n-2)t/2}u_io\exp_{y_i}(e^t\theta), \,\, {\rm et} \,\, a(y_i,t,\theta)=\log J(y_i,e^t,\theta). $$ 

\underbar {\bf Lemma 1:}

\bigskip

The function $ w_i $ is solution of:

$$  -\partial_{tt} w_i-\partial_t a \partial_t w_i-L_{\theta}(y_i,e^t)+c w_i=n(n-2)w_i^{N-1}, $$
 
with,

 $$ c = c(y_i,t,\theta)=\left ( \dfrac{n-2}{2} \right )^2+ \dfrac{n-2}{2} \partial_t a  -\lambda e^{2t}, $$

\underbar {\bf Proof:}

\smallskip

We write:

$$ \partial_t w_i=e^{nt/2}\partial_r \bar u_i+\dfrac{n-2}{2} w_i,\,\, \partial_{tt} w_i=e^{(n+2)t/2} \left [\partial_{rr} \bar u_i+\dfrac{n-1}{e^t}\partial_r \bar u_i \right ]+\left ( \dfrac{n-2}{2} \right )^2 w_i. $$

$$ \partial_t a =e^t\partial_r \log J(y_i,e^t,\theta), \partial_t a \partial_t w_i=e^{(n+2)t/2}\left [ \partial_r \log J\partial_r \bar u_i \right ]+\dfrac{n-2}{2} \partial_t a w_i.$$

the lemma  is proved.

\bigskip

Now we have, $ \partial_t a=\dfrac{ \partial_t b_1}{b_1} $, $ b_1(y_i,t,\theta)=J(y_i,e^t,\theta)>0 $,

\bigskip

We can write,

$$ -\dfrac{1}{\sqrt {b_1}}\partial_{tt} (\sqrt { b_1} w_i)-L_{\theta}(y_i,e^t)w_i+[c(t)+ b_1^{-1/2} b_2(t,\theta)]w_i=n(n-2){w_i}^{N-1}, $$

where, $ b_2(t,\theta)=\partial_{tt} (\sqrt {b_1})=\dfrac{1}{2 \sqrt { b_1}}\partial_{tt}b_1-\dfrac{1}{4(b_1)^{3/2}}(\partial_t b_1)^2 .$

\bigskip

Let,

$$ \tilde w_i=\sqrt {b_1} w_i, $$

\underbar {\bf Lemma 2:}

\smallskip

The function $ \tilde w_i $ is solution of:

$$ -\partial_{tt} \tilde w_i+\Delta_{g_{y_i, e^t, {}_{{\mathbb S}_{n-1}}}} (\tilde w_i)+2\nabla_{\theta}(\tilde  w_i) .\nabla_{\theta} \log (\sqrt {b_1})+(c+b_1^{-1/2} b_2-c_2) \tilde w_i= $$

$$ = n(n-2)\left (\dfrac{1}{b_1} \right )^{(N-2)/2} {\tilde w_i}^{N-1}, $$

where, $ c_2 =[\dfrac{1}{\sqrt {b_1}} \Delta_{g_{y_i, e^t, {}_{{\mathbb S}_{n-1}}}}(\sqrt{b_1}) + |\nabla_{\theta} \log (\sqrt {b_1})|^2] . $

\bigskip

\underbar {\bf Proof:}

\bigskip

We have: 

$$ -\partial_{tt} \tilde w_i-\sqrt {b_1} \Delta_{g_{y_i, e^t, {}_{{\mathbb S}_{n-1}}}} w_i+(c+b_2) \tilde w_i= n(n-2)\left (\dfrac{1}{b_1} \right )^{(N-2)/2} {\tilde w_i}^{N-1}, $$

But,

$$ \Delta_{g_{y_i, e^t, {}_{{\mathbb S}_{n-1}}}} (\sqrt {b_1} w_i)=\sqrt {b_1} \Delta_{g_{y_i, e^t, {}_{{\mathbb S}_{n-1}}}} w_i-2 \nabla_{\theta} w_i .\nabla_{\theta} \sqrt {b_1}+ w_i \Delta_{g_{y_i, e^t, {}_{{\mathbb S}_{n-1}}}}(\sqrt {b_1}), $$

and,

$$ \nabla_{\theta} (\sqrt {b_1} w_i)=w_i \nabla_{\theta} \sqrt {b_1}+ \sqrt {b_1} \nabla_{\theta} w_i, $$

we deduce than,

$$  \sqrt {b_1} \Delta_{g_{y_i, e^t, {}_{{\mathbb S}_{n-1}}}} w_i= \Delta_{g_{y_i, e^t, {}_{{\mathbb S}_{n-1}}}} (\tilde w_i)+2\nabla_{\theta}(\tilde  w_i) .\nabla_{\theta} \log (\sqrt {b_1})-c_2 \tilde w_i, $$

with $ c_2=[\dfrac{1}{\sqrt {b_1}} \Delta_{g_{y_i, e^t, {}_{{\mathbb S}_{n-1}}}}(\sqrt{b_1}) + |\nabla_{\theta} \log (\sqrt {b_1})|^2] . $ The lemma is proved.

\bigskip

\underbar {\bf The "moving-plane" method:}

\bigskip

Let $ \xi_i $ a real number,  and suppose $ \xi_i \leq t $. We set $ t^{\xi_i}=2\xi_i-t $ and $ \tilde w_i^{\xi_i}(t,\theta)=\tilde w_i(t^{\xi_i},\theta) $.

\bigskip

We have, 

$$ -\partial_{tt} \tilde w_i^{\xi_i}+\Delta_{g_{y_i, e^{t^{\xi_i}} {}_{{\mathbb S}_{n-1}}}} (\tilde w_i)+2\nabla_{\theta}(\tilde  w_i^{\xi_i}) .\nabla_{\theta} \log (\sqrt {b_1}) \tilde w_i^{\xi_i}+[c(t^{\xi_i})+b_1^{-1/2}(t^{\xi_i},.)b_2(t^{\xi_i})-c_2^{\xi_i}] \tilde w_i^{\xi_i}= $$

$$ =n(n-2)\left (\dfrac{1}{b_1^{\xi_i}} \right )^{(N-2)/2} { ({\tilde w_i}^{\xi_i}) }^{N-1}. $$

By using the same arguments than in [B2], we have:

\bigskip

\underbar {\bf Proposition 3:}

\smallskip

We have:

$$ 1)\,\,\, \tilde w_i(\lambda_i,\theta)-\tilde w_i(\lambda_i+4,\theta) \geq \tilde k>0, \,\, \forall \,\, \theta \in {\mathbb S}_{n-1}. $$

For all $ \beta >0 $, there exist $ c_{\beta} >0 $ such that:

$$ 2) \,\,\, \dfrac{1}{c_{\beta}} e^{(n-2)t/2}\leq \tilde w_i(\lambda_i+t,\theta) \leq c_{\beta}e^{(n-2)t/2}, \,\, \forall \,\, t\leq \beta, \,\, \forall \,\, \theta \in {\mathbb S}_{n-1}. $$

We set,

$$ \bar Z_i=-\partial_{tt} (...)+\Delta_{g_{y_i, e^t, {}_{{\mathbb S}_{n-1}}}} (...)+2\nabla_{\theta}(...) .\nabla_{\theta} \log (\sqrt {b_1})+(c+b_1^{-1/2} b_2-c_2)(...) $$

{\bf Remark:} In the operator $ \bar Z_i $, by using the proposition 3, the coeficient $ c+b_1^{-1/2}b_2-c_2 $ satisfies:

$$ c+b_1^{-1/2}b_2-c_2 \geq k'>0,\,\, {\rm pour }\,\, t<<0, $$

it is fundamental if we want to apply Hopf maximum principle.

\bigskip

\underbar {\bf Goal:}

\bigskip

Like in [B2], we have elliptic second order operator. Here it is $ \bar Z_i $, the goal is to use the "moving-plane" method to have a contradiction. For this, we must have:

$$ \bar Z_i(\tilde w_i^{\xi_i}-\tilde w_i) \leq 0, \,\, {\rm if} \,\, \tilde w_i^{\xi_i}-\tilde w_i \leq 0. $$

We write:

$$ \bar Z_i(\tilde w_i^{\xi_i}-\tilde w_i)= (\Delta_{g_{y_i, e^{t^{\xi_i}}, {}_{{\mathbb S}_{n-1}}}}-\Delta_{g_{y_i, e^{t}, {}_{{\mathbb S}_{n-1}}}}) (\tilde w_i^{\xi_i})+ $$

$$ +2(\nabla_{\theta, e^{t^{\xi_i}}}-\nabla_{\theta, e^t})(w_i^{\xi_i}) .\nabla_{\theta, e^{t^{\xi_i}}} \log (\sqrt {b_1^{\xi_i}})+ 2\nabla_{\theta,e^t}(\tilde w_i^{\xi_i}).\nabla_{\theta, e^{t^{\xi_i}}}[\log (\sqrt {b_1^{\xi_i}})-\log \sqrt {b_1}]+ $$ 

$$ +2\nabla_{\theta,e^t} w_i^{\xi_i}.(\nabla_{\theta,e^{t^{\xi_i}}}-\nabla_{\theta,e^t})\log \sqrt {b_1}- [(c+b_1^{-1/2} b_2-c_2)^{\xi_i}-(c+b_1^{-1/2}b_2-c_2)]\tilde w_i^{\xi_i} + $$

$$ + n(n-2)\left ( \dfrac{1}{b_1^{\xi_i}} \right )^{(N-2)/2} ({\tilde w_i}^{\xi_i})^{N-1}-n(n-2)\left ( \dfrac{1}{b_1} \right )^{(N-2)/2} {\tilde w_i}^{N-1}.\qquad (***1) $$

Clearly, we have:

\bigskip

\underbar {\bf Lemma 3 :}

$$ b_1(y_i,t,\theta)=1-\dfrac{1}{3} Ricci_{y_i}(\theta,\theta)e^{2t}+\ldots, $$

$$ R_g(e^t\theta)=R_g(y_i) + <\nabla R_g(y_i)|\theta > e^t+\dots . $$

According to proposition 1 and lemma 3,

\bigskip

\underbar {\bf Propostion 4 :}

$$ \bar Z_i(\tilde w_i^{\xi_i}-\tilde w_i) \leq {b_1}^{(2-N)/2}[(\tilde w_i^{\xi_i})^{N-1}- \tilde w_i^{N-1}]+  $$

$$ +C|e^{2t}-e^{2t^{\xi_i}}|\left [|\nabla_{\theta} {\tilde w_i}^{\xi_i}| + |\nabla_{\theta}^2(\tilde w_i^{\xi_i})|+ |Ricci_{y_i}|[\tilde w_i^{\xi_i}+(\tilde w_i^{\xi_i})^{N-1}] + |R_g(y_i)| \tilde w_i^{\xi_i} \right ] + C'|e^{3t^{\xi_i}}-e^{3t}|. $$

\underbar {\bf Proof:}

\bigskip

We use proposition 1, we have:

$$ a(y_i,t,\theta)=\log J(y_i,e^t,\theta)=\log b_1, |\partial_t b_1(t)|+|\partial_{tt} b_1(t)|+|\partial_{tt} a(t)|\leq C e^{2t}, $$

and,

$$ |\partial_{\theta_j} b_1|+|\partial_{\theta_j,\theta_k} b_1|+\partial_{t,\theta_j}b_1|+|\partial_{t,\theta_j,\theta_k} b_1|\leq C e^{2t}, $$

then,

$$ |\partial_t b_1(t^{\xi_i})-\partial_t b_1(t)|\leq C'|e^{2t}-e^{2t^{\xi_i}}|,\,\, {\rm on} \,\, ]-\infty, \log \epsilon_1]\times {\mathbb S}_{n-1},\forall \,\, x\in B(x_0,\epsilon_1) $$

Locally,

$$ \Delta_{g_{y_i, e^t, {}_{{\mathbb S}_{n-1}}}}= L_{\theta}(y_i,e^t)=-\dfrac{1}{\sqrt {|{\tilde g}^k(e^t,\theta)|}}\partial_{\theta^l}[{\tilde g}^{\theta^l \theta^j}(e^t,\theta)\sqrt { |{\tilde g}^k(e^t,\theta)|}\partial_{\theta^j}] . $$

Thus, in $ [0,\epsilon_1]\times U^k $, we have,

$$ A_i=\left [{ \left [ \dfrac{1}{\sqrt {|{\tilde g}^k|}}\partial_{\theta^l}({\tilde g}^{\theta^l \theta^j}\sqrt { |{\tilde g}^k|}\partial_{\theta^j}) \right ] }^{\xi_i}- \dfrac{1}{\sqrt {|{\tilde g}^k|}}\partial_{\theta^l}({\tilde g}^{\theta^l \theta^j}\sqrt { |{\tilde g}^k|}\partial_{\theta^j}) \right ](\tilde w_i^{\xi_i}) $$

then, $ A_i=B_i+D_i $ with,

$$ B_i=\left [ {\tilde g}^{\theta^l \theta^j}(e^{t^{\xi_i}}, \theta)-{\tilde g}^{\theta^l \theta^j}(e^t,\theta) \right ] \partial_{\theta^l \theta^j}\tilde w_i^{\xi_i}(t,\theta), $$

and,

$$ D_i=\left [ \dfrac{1}{ \sqrt {| {\tilde g}^k|}(e^{t^{\xi_i}},\theta )  }           \partial_{\theta^l}[{\tilde g }^{\theta^l \theta^j}(e^{t^{\xi_i}},\theta)\sqrt {| {\tilde g}^k|}(e^{t^{\xi_i}},\theta)  ] -\dfrac{1}{ \sqrt {| {\tilde g}^k|}(e^t,\theta) } \partial_{\theta^l} [{\tilde g }^{\theta^l \theta^j}(e^t,\theta)\sqrt {| {\tilde g}^k|}(e^t,\theta) ] \right ] \partial_{\theta^j} \tilde w_i^{\xi_i}(t,\theta), $$

we deduce,

$$ A_i \leq C_k|e^{2t}-e^{2t^{\xi_i}}|\left [ |\nabla_{\theta} \tilde w_i^{\xi_i}| + |\nabla_{\theta}^2(\tilde w_i^{\xi_i})| \right ], $$

If we take $ C=\max \{ C_i, 1 \leq i\leq q \} $ and if w use $ (***1) $, we obtain proposition 4.

We have,

$$ c(y_i,t,\theta)=\left ( \dfrac{n-2}{2} \right )^2+ \dfrac{n-2}{2} \partial_t a + R_g e^{2t}, \qquad (\alpha_1) $$ 

$$ b_2(t,\theta)=\partial_{tt} (\sqrt {b_1})=\dfrac{1}{2 \sqrt { b_1}}\partial_{tt}b_1-\dfrac{1}{4(b_1)^{3/2}}(\partial_t b_1)^2 ,\qquad (\alpha_2) $$ 

$$ c_2=[\dfrac{1}{\sqrt {b_1}} \Delta_{g_{y_i, e^t, {}_{{\mathbb S}_{n-1}}}}(\sqrt{b_1}) + |\nabla_{\theta} \log (\sqrt {b_1})|^2], \qquad (\alpha_3) $$

Then,

$$ \partial_{t}c(y_i,t,\theta)=\dfrac{(n-2)}{2}\partial_{tt}a+2e^{2t}R_g(e^t\theta)+e^{3t}<\nabla R_g(e^t\theta)|\theta >, $$

by proposition 1,

$$ |\partial_tc_2|+|\partial_t b_1|+|\partial_t b_2|+|\partial_t c|\leq K_1e^{2t}, $$

\underbar {\bf The case: $ 0 <m \leq \lambda+R_g \leq \dfrac{1}{m} $ for the equation $ \Delta u-\lambda u=n(n-2)u^{N-1} $ }

\bigskip

Let $ x_0 $ a point of $ M $, we consider a conformal change of metric $ \tilde g= \phi^{4/(n-2)} g $ such that, $ \tilde Ricci(x_0)=0 $. See for example [Au] ( also Lee and Parker [L,P]).

\bigskip

We are concerning by the following equation,

$$ \Delta_g u-\lambda u=n(n-2)u^{N-1}, $$

the conformal change of metric give when we set $ v=u/\phi $,

$$ \Delta_{\tilde g} v +\tilde R_{\tilde g} v =n(n-2)v^{N-1}+(\lambda + \tilde R_g) \phi^{N-2} v. $$

The notation $ \tilde R $ is for $ \dfrac{n-2}{4(n-1)} R $ and $ R=R_g $ or $ R=R_{\tilde g} $.

\bigskip

Our calculus for the metric $ \tilde g $ are the same that for the metric $ g $. But we have some new properties:

$$ \sqrt { det (\tilde g_{y_i,jk})}=1-\dfrac{1}{3}\tilde Ricci(y_i)(\theta,\theta)r^2+..., \,\, {\rm and} \,\, \tilde R_{\tilde g}(y_i) \to 0, \,\, \tilde Ricci(y_i) \to 0. $$

If we see the coeficient in the term $ e^{2t^{\xi_i}}-e^{2t} $, we can say that all those  terms are tending to 0, see proposition 4. Only the term $ (\lambda+ \tilde R_g)(e^{2t^{\xi_i}}-e^{2t}) \leq m (e^{2t^{\xi}}-e^{2t}) $ ( $ m>0 $), is the biggest.

\bigskip

In fact, the increment of the local expression of the metric $ {\tilde g}_{jk}^{\xi_i}-{\tilde g}_{jk} $, have terms of type $ \partial_{\theta_j} \tilde w_i^{\xi_i} $ et $ \partial_{\theta_j,\theta_k} \tilde w_i^{\xi_i} $ but we know by proposition 2 that those terms tend to 0 because the limit function is radial and do not depend on the angles.

\bigskip

We apply  proposition 3. We take $ t_i=\log \sqrt {l_i} $ with $ l_i $ like in proposition 2. The fact $ \sqrt {l_i} [u_i(y_i)]^{2/(n-2)} \to +\infty $ ( see proposition 2), implies $ t_i=\log \sqrt {l_i} > \dfrac{2}{n-2} \log u_i(y_i) + 2 =\lambda_i+2 $. Finaly, we can work  on $ ]-\infty, t_i] $.

\bigskip

We define $ \xi_i $ by:

$$ \xi_i=\sup \{ \lambda \leq \lambda_i+2, \,\, \tilde w_i(2\lambda-t,\theta)-\tilde w_i(t,\theta)\leq 0\,\, {\rm on} \,\, [\lambda,t_i]\times {\mathbb S}_{n-1} \}. $$ 

If we use proposition 4 and the similar technics that in [B2] we can deduce by Hopf maximum principle,

$$ \max_{{\mathbb S}_{n-1}} \tilde w_i(t_i,\theta) \leq \min_{{\mathbb S}_{n-1}} \tilde w_i(2\xi_i-t_i,\theta), $$

which implies,

$$ {l_i}^{(n-2)/2} u_i(y_i) \times \min_M u_i \leq c. $$

It is in contradiction with  proposition 2.

\bigskip

Then we have,

$$ \sup_K u \times \inf_M u \leq c=c(K,M,m,g,n). $$

\underbar {\bf Application:}

\bigskip

Let $ M $ a Riemannian manifold of dimension $ n\geq 3 $, and consider a sequence of functions $ u_i $ such that:

$$ \Delta u_i+\epsilon_i u_i=n(n-2){u_i}^{N-1}, \,\, \epsilon_i \to 0 $$

If, the scalar curvature $ R_g \geq m >0 $ on $ M $, then, applying the previous result with $ \lambda=-\epsilon_i $, we obtain:

$$ \sup_M u_i \times \inf_M u_i \leq c, \,\, \forall \,\, i, $$

\underbar {\it Proof of the theorem 4:}

Without loss of generality we suppose,

$$ \Delta u_i+\epsilon_i u_i=u_i^{N-1}, \,\, {\rm et} \,\, \max_M u_i \to 0. $$

\underbar {\bf Lemma 1:} There exist a positive constant, $ c $ such that:

$$ \sup_M u_i \leq c \inf_M u_i, \,\, \forall \, i. $$

\underbar {\bf Proof of lemma 1:}

\bigskip

Suppose by contradiction:

$$ \limsup_{i \to +\infty} \dfrac{\sup_M u_i}{\inf_M u_i} =+\infty, $$

After passing to a subsequence, we can assume: $ \dfrac{\sup_M u_i}{\inf_M u_i} \to +\infty. $

\bigskip

We have, $ \sup_M u_i=u_i(y_i) $ et $ \inf_M u_i=u_i(x_i) $. We also suppose, $ x_i \to x $ et $ y_i \to y $.

\bigskip

Let $ L $ be a minimizing curve between $ x $ and $ y $, take $ \delta > 0 $ such that $ \delta < inj_g(M) $, with $ inj_g(M) $ the injectivity radius of the compact manifold $ M $.

\bigskip

For all $ a \in L $, $ [B(a,\delta), (\exp_a)^{-1}] $ is a local chart around  $ a $, but $ L $ is compact. We can cover this curve by a finite number of balls centred in a points of $ L $ and of radius $ \delta / 5 $. Let $ a_1,\ldots, a_k $ those points, with, $ a_1=x $ and $ a_k=y $.

\bigskip

In each ball $ B(a_j,\delta) $, $ u_i $ is solution of, $ \Delta u_i+(\epsilon_i-u_i^{N-2})u_i=0 $, we use the fact $ \sup_M u_i \to 0 $ and we apply the Harnack inequality of [G-T] ( see theorem 8.20), we obtain:

$$ \sup_{B(a_j,\delta/5)} u_i \leq C_j \inf_{B(a_j,\delta/5)} u_i, \,\, j=1,\ldots,k. $$

We deduce:

$$ \sup_{B(y,\delta/5)} u_i \leq C_kC_{k-1}.\ldots. C_1\inf_{B(x,\delta/5)} u_i, $$

In other words:

$$ \sup_M u_i \leq C_k.\dots.C_1 \inf_M u_i. $$

It's in contradiction with our hypothesis.

\bigskip

\underbar {\bf Lemma 2:} There exist two constants, $ k_1,k_2>0 $ such that:

$$ k_1 {\epsilon_i}^{(n-2)/4} \leq u_i(x) \leq k_2 {\epsilon_i}^{(n-2)/4}, \,\, \forall \,\, x\in M, \,\, \forall \,\, i. $$

\underbar {\bf Proof of lemma 2:}

\bigskip

Let $ G_i $ the Green function of the operator $ \Delta+\epsilon_i $, this equation satisfies:

$$ \int_M G_i(x,y)dV_g(y)=\dfrac{1}{\epsilon_i}, \,\, \forall \,\, x\in M. $$

We write:

$$ \inf_M u_i=u_i(x_i)=\int_M G_i(x_i,y)u_i^{N-1}(y)dV_g(y) \geq (\inf_M u_i)^{N-1}\int_M G_i(x_i,y)dV_g(y) =\dfrac{(\inf_M u_i)^{N-1}}{\epsilon_i}, $$

thus,

$$ \inf_M u_i \leq {\epsilon_i}^{(n-2)/4}. $$

We the same idea we can prove, $ \sup_M u_i \geq {\epsilon_i}^{(n-2)/4} $. We deduce lemma 2 from lemma 1 and the two last inequalities.

\bigskip

\underbar {\bf Lemma 3:} There exist a rank $ i_0 $ such that, $ u_i \equiv {\epsilon_i}^{(n-2)/4}. $ for $ i \geq i_0 $.

\bigskip

\underbar{ \bf Proof of lemma 3:}

\bigskip

Let, $ w_i=\dfrac{u_i}{{\epsilon_i}^{(n-2)/4}} $. This function is solution of:

$$ \Delta w_i = \epsilon_i(w_i^{N-1}-w_i)=\epsilon_i w_i(w_i^{N-2}-1). \qquad (*) $$

\underbar {\bf Case 1: $  N-2 \geq 1 $ ($ 3\leq n \leq 6 $),}

\smallskip

To simplify our computations we suppose that $ N-2 $ is an integer.

\bigskip

According to binomial formula, $ w_i^{N-2}-1=(w_i-1)(1+w_i+...) $, we multiply $ (*) $ by $ w_i-1 $ and we integrate, we obtain:

$$ \int_M |\nabla w_i|^2 \leq C \epsilon_i \int_M |w_i-1|^2, $$

Suppose that we have infinity $ i $, such that $ w_i \not \equiv 1 $, then we can consider the following functions: $ z_i=\dfrac{w_i-1}{||w_i-1||_2} $.

\bigskip

$ z_i $ verifiy, $ ||z_i||_2=1, ||\nabla z_i||_2^2\leq C\epsilon_i \to 0 $, thus, $ z_i \to 1 $ in $ L^2(M) $ and in particular, $ \int_M z_i w_i(1+w_i+...) \to C'\not = 0 $ ( by using lemma 2). But, if we integrate $ (*) $, we find $ \int_M z_i w_i(1+w_i+....)=0 $, it's a contradiction.

\bigskip

Thus, there exist a rank such that $ w_i \equiv 1 $ after this rank.

\bigskip

\underbar {\bf Case 2: $ 0 < N-2 < 1 $ ( $ n\geq 7 $):}

\smallskip

To simplify our computations, we suppose that $ 1/(N-2) $ is an integer.

\bigskip

Now we take $ w_i^{N-2}-1 $ and we write $ w_i-1=(w_i^{N-2})^{1/(N-2)}-1 $, by using the binomial formula and the same ideas than in the previous case we obtain our result.

\bigskip

\bigskip

\bigskip

\bigskip

\begin{center}

{\bf  ACKNOWLEDGEMENT. }

\end{center}

\smallskip

This work was done when the author was in Greece at Patras. The author is grateful to Professor Athanase Cotsiolis, the Department of Mathematics of Patras University  and the IKY Foundation for hospitality and the very good conditions of work. 

\bigskip

\underbar {\bf References:}

\bigskip

[Au] T. Aubin. Some Nonlinear Problems in Riemannian Geometry. Springer-Verlag 1998.

\smallskip

[B1] S.S Bahoura, Diff\'erentes Estimations du $ \sup u \times \inf u $ pour l'\'Equation de la Courbure Sclaire Prescrite en dimension $ n\geq 3 $. J. Math. Pures . Appl. (9)82 (2003), no.1, 43-66.

\smallskip

[B2] S.S Bahoura. Majorations du type $ \sup u \times \inf u \leq c $ pour l'\'equation de la courbure scalaire sur un ouvert de $ {\mathbb R}^n, n\geq 3 $. J. Math. Pures. Appl.(9) 83 2004 no, 9, 1109-1150.

\smallskip

[B-L-S] H. Brezis, Yy. Li Y-Y, I. Shafrir. A sup+inf inequality for some
nonlinear elliptic equations involving exponential
nonlinearities. J.Funct.Anal.115 (1993) 344-358.

\smallskip

[B-M] H.Brezis and F.Merle, Uniform estimates and blow-up bihavior for solutions of $ -\Delta u=Ve^u $ in two dimensions, Commun Partial Differential Equations 16 (1991), 1223-1253.

\smallskip

[B-V, V] M-F. Bidaut-V\'eron, L. V\'eron. Nonlinear elliptic equations on compact Riemannian manifolds and asymptotics od Emden equations. Invent.Math. 106 (1991), no3, 489-539.

\smallskip

[C-G-S], L. Caffarelli, B. Gidas, J. Spruck. Asymptotic symmetry and local
behavior of semilinear elliptic equations with critical Sobolev
growth. Comm. Pure Appl. Math. 37 (1984) 369-402.

\smallskip

[C-L 1] A sharp sup+inf inequality for a nonlinear elliptic equation in ${\mathbb R}^2$.
Commun. Anal. Geom. 6, No.1, 1-19 (1998).

\smallskip

[C-L 2] C-C.Chen, C-S. Lin. Estimates of the conformal scalar curvature
equation via the method of moving planes. Comm. Pure
Appl. Math. L(1997) 0971-1017.

\smallskip

[D-H-R] O. Druet, E. Hebey, F.Robert, Blow-up theory in Riemannian Geometry, Princeton University Press.

\smallskip

[He] E. Hebey, Analyse non lineaire sur les Vari\'et\'e, Editions Diderot.

\smallskip

[He,V] The best constant problem in the Sobolev embedding theorem for complete Riemannian manifolds.  Duke Math. J.  79  (1995),  no. 1, 235--279.

\smallskip

[L,P] J.M. Lee, T.H. Parker. The Yamabe problem. Bull.Amer.Math.Soc (N.S) 17 (1987), no.1, 37 -91.

\smallskip

[L 1] Yy. Li. Harnack Type Inequality: the Method of Moving Planes. Commun. Math. Phys. 200,421-444 (1999).

\smallskip

[L 2] Yy. Li. Prescribing scalar curvature on $ {\mathbb S}_n $ and related
Problems. C.R. Acad. Sci. Paris 317 (1993) 159-164. Part
I: J. Differ. Equations 120 (1995) 319-410. Part II: Existence and
compactness. Comm. Pure Appl.Math.49 (1996) 541-597.

\smallskip

[L-Z] Yy. Li, L. Zhang. A Harnack type inequality for the Yamabe equation in low dimensions.  Calc. Var. Partial Differential Equations  20  (2004),  no. 2, 133--151.

\smallskip

[S] I. Shafrir. A sup+inf inequality for the equation $ -\Delta u=Ve^u $. C. R. Acad.Sci. Paris S\'er. I Math. 315 (1992), no. 2, 159-164.

\end{document}